\documentclass[11pt]{amsart}
\usepackage{amssymb,amsmath}
\usepackage{subfigure}
\usepackage{graphpap,latexsym,epsf}
\usepackage{color,graphicx,psfrag}

\def\supp(#1){\{#1>0\}}
\def\k{\varkappa}
\def\vk{\varkappa}

\newcommand{\alchi}{\raisebox{1.7pt}{$\chi$}}
\newcommand{\R}{{\mathbb R}}
\newcommand{\N}{{\mathbb N}}

\newcommand{\e }{\varepsilon }
\newcommand{\eps }{\varepsilon }

\newcommand{\sgn}{\mathop{\rm sgn}}

\newenvironment{pf}{\noindent{\sc Proof}.\enspace}{\hfill\qed
\medskip}
\newenvironment{pfn}[1]{\noindent{\sc Proof of {#1}.\enspace}}{\hfill\qed
\medskip}
\newtheorem{Theorem}{Theorem}[section]

\newtheorem{Lemma}[Theorem]{Lemma}
\newtheorem{Proposition}[Theorem]{Proposition}

\newtheorem{remark}[Theorem]{Remark}
\theoremstyle{definition}

\begin{document}

\title[Minimal coexistence configurations for multispecies
systems]{Minimal coexistence configurations for multispecies systems}

\author[Monica Conti]{Monica Conti}
\address{\hbox{\parbox{5.7in}{\medskip\noindent{Monica Conti:
        Politecnico di Milano, Dipartimento di Matematica
        ``F. Brioschi'', Via Bonardi 9, I-20133 Milano, Italy.
        \em{E-mail address: }{\tt monica.conti@polimi.it.}}}}}
\author[Veronica Felli]{Veronica Felli}
\address{\hbox{\parbox{5.7in}{\medskip\noindent{Veronica Felli:
        Universit\`a di Milano Bicocca, Dipartimento di
        Ma\-t\-ema\-ti\-ca e Applicazioni, Via Cozzi 53, 20125 Milano,
        Italy. \em{E-mail address: }{\tt veronica.felli@unimib.it}.}}}}

\date{November 18, 2008}

\thanks{Supported by Italy MIUR, national project ``Variational
  Methods and Nonlinear Differential Equations''.
  \\
  \indent 2000 {\it Mathematics Subject Classification.} 35J55, 47J30,
  92D25.\\
  \indent {\it Keywords.} Segregation states, competing multispecies
  systems, domain perturbation.}

 \begin{abstract}
   \noindent We deal with strongly competing multispecies systems of
   Lotka-Volterra type with homogeneous Neumann boundary conditions in
   dumbbell-like domains. Under suitable non-degeneracy assumptions,
   we show that, as the competition rate grows indefinitely, the
   system reaches a state of coexistence of all the species in spatial
   segregation. Furthermore, the limit configuration is a local
   minimizer for the associated free energy.
 \end{abstract}

 \maketitle

\section{Introduction}\label{sec:assumpt-main-results}
In this paper we consider  the system of
$k\geq 2$ elliptic equations
\begin{equation}\label{modelLK}
-\Delta u_i+u_i= f_i(u_i)-\vk \;u_i\sum_{j\neq i}u_j^2,\quad\text{ in } \Omega,
\end{equation}
for $i=1,\dots,k$.  It models the steady states of $k$ organisms, each
of density $u_i$, which coexist in a smooth, connected, bounded domain
$\Omega\subset\R^N$; their dynamics is ruled out by internal growth $f_i$'s and 
mutual competition of Lotka-Volterra type with parameter $\vk>0$.
Systems of this form have attracted considerable attention both in
ecology and social science since they furnish a relatively simple
model to study the behavior of $k$ populations competing for the same
resource $\Omega$.  One of the main question is to investigate whether
\emph{coexistence} may occur, namely the existence of equilibrium
configurations where all the densities $u_i$ are strictly positive on
sets of positive measure, or the internal dynamic leads to
\emph{extinction}, that is steady states where one or more densities
are null.  Many results are nowadays available, dealing mainly with
$k=2$ populations.  We quote among others
\cite{eflg,gl,kl,lm,mhmm,skt}, where for logistic internal growth
$f_i(u)=u(a_i-u)$, both the situation are proved to be possible
depending on the relations between the diffusion rates and the
coefficients of intra--specific and
of inter--specific competitions, see also \cite{ddh,dd}.

A different perspective is proposed in
\cite{cf,ctv2,ctv-asymp,dd3,dg1,dhmp}, where the authors study the
effect of very strong competition, letting the parameter $\vk$
growing indefinitely.  It is observed (see
Section~\ref{sec:asymptotic-analysis}) that the presence of large
interactions of competitive type produces, in the limit configuration
as $\vk\to \infty$, the \emph{spatial segregation} of the
densities, meaning that if $(u_i^{\vk})_{i=1,\dots,k}$
solves~(\ref{modelLK}), then $u_i^{\vk}$ converges (in a suitable
sense) to some $u_i$ which satisfies
\begin{equation}\label{segrego}
  u_i(x)\cdot u_j(x)=0\text{ a.e. in }\Omega,\qquad\text{for all } i\neq j.
\end{equation}
A number of qualitative properties of the possible coexistence states
$u_i$ and their supports is proved in
\cite{ctv-var,ctv-asymp,ctv-uniq}, with the aim of describing the way
the territory is partitioned by the segregated populations. We refer
the interested reader to the above quoted papers for  details on the
regularity theory so far developed and to \cite{ctv2,ctv-fucik} for
some applications.

A further point of interest is to establish if coexistence of the
species is possible in a segregated configuration: do all the species
survive when the intra specific competition becomes larger and larger?
The answer cannot be positive in general: \cite{kw} shows that in any
\emph{convex domain} the only stable configurations are those where
only one specie is alive.  It is worth pointing out that in
\cite{ctv-var,ctv-asymp}, the strict positivity of each component in
the limiting configuration is guaranteed by simply forcing
non-homogeneous Dirichlet boundary conditions
\begin{equation}\label{dirnonomo}
u_i=\phi_i\quad \text{on }\partial\Omega,
\end{equation}
with $\phi_i>0$ on a set of positive $(N-1)$--measure.  Coexistence
results for competing systems under more natural homogeneous boundary
conditions are obtained in \cite{cf} for the Dirichlet case
\begin{equation}\label{diromo}
u_i=0\quad\text{on }\partial\Omega,
\end{equation}
with interactions of the form $\vk u_i\sum_{j\neq i}u_j$.  To avoid the extinction
predicted by \cite{kw}, a special class of
non-convex domains close to a union of $k$ disjoint balls is considered. Suitable
non-degeneracy assumptions on the $f_i$'s allow the application of a
domain perturbation technique envisaged in \cite{d} which strongly relies on the continuity of
the eigenvalues of the Laplace operator with respect to the domain.
It is well known that such a property does not hold in the case of Neumann boundary conditions,
see for instance \cite{arrieta}.  Hence, in order to treat Neumann
no-flux boundary conditions, a different approach is needed.

\begin{figure}[h]
  \centering \subfigure [example of $\Omega_{0}$ with $k=3$]{
    \begin{psfrags}
      \psfrag{B1}{$\Omega^1$} \psfrag{B2}{$\Omega^2$}
      \psfrag{B3}{$\Omega^3$}
      \includegraphics[width=5.5cm]{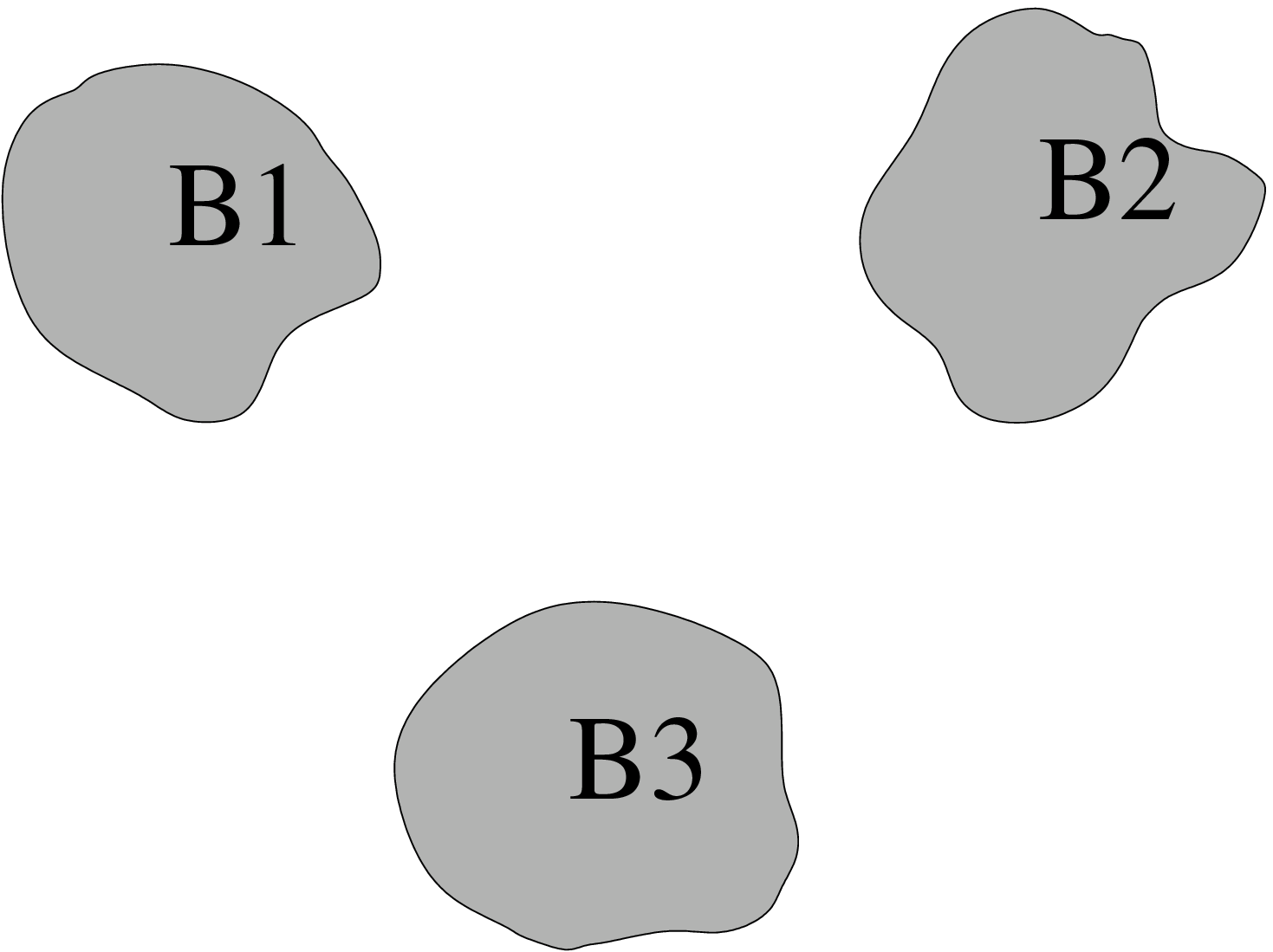}
    \end{psfrags}
  } \quad\quad\quad \subfigure[example of $\Omega_\e$ approximating $\Omega_0$ ] {
    \includegraphics[width=5.5cm]{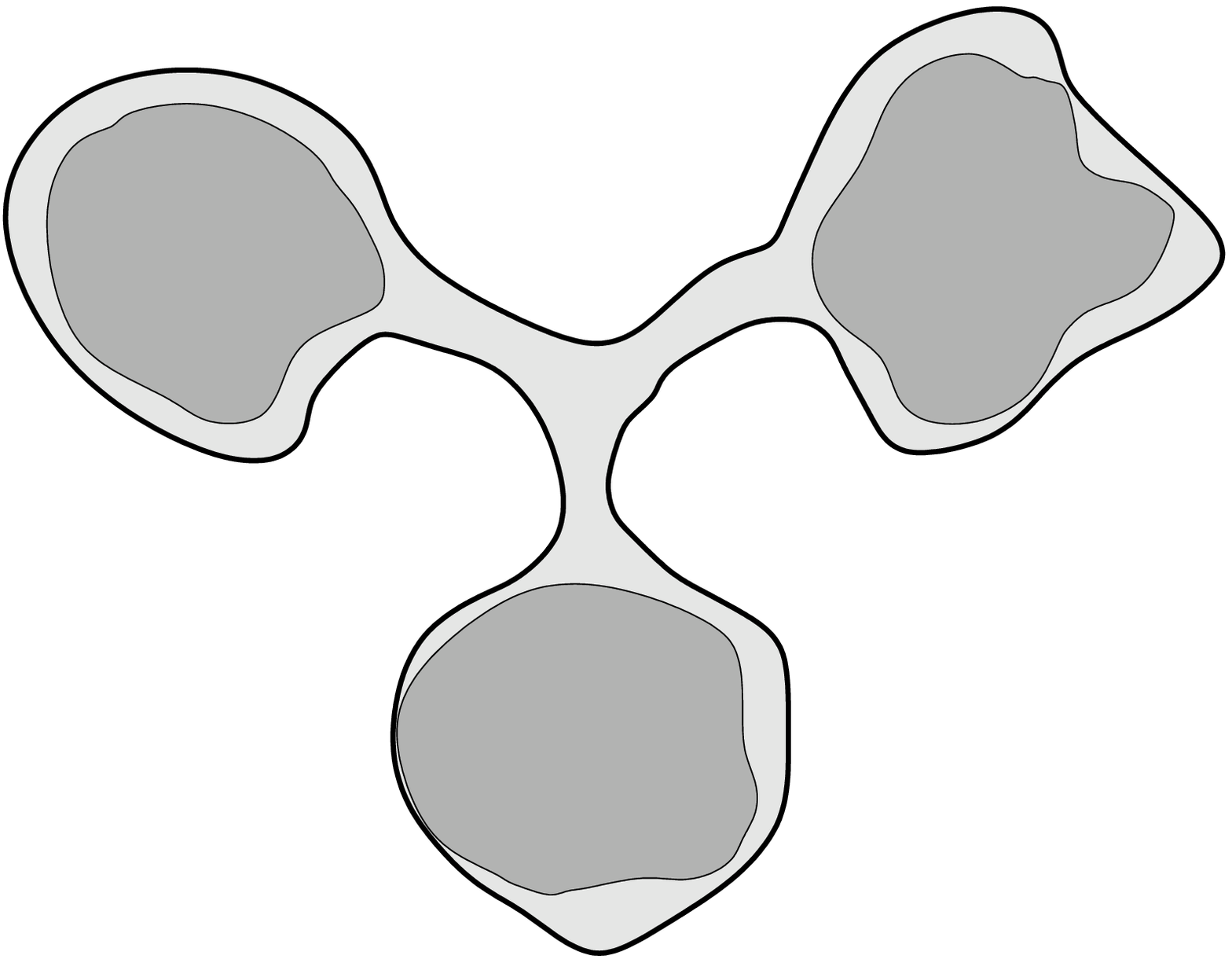}
  }
  \caption{}\label{fig:dd}
\end{figure}

This is precisely the aim of the present paper: we deal with system
\eqref{modelLK} coupled with
\begin{equation}\label{neuomo}
\frac{\partial u_i}{\partial \nu}=0\quad\text{on }\partial\Omega,
\end{equation}
in
a class of non-convex domains $\Omega=\Omega_\e$ suitably approximating a given domain $\Omega_0$ 
composed by $k$ disjoint open sets, see Figure
\ref{fig:dd}.

Due to the variational structure of  problem \eqref{modelLK},
the following  \emph{free energy} functional
\begin{equation}\label{freeenergy}
  J_{\Omega}(U)=\sum_{i=1,\dots, k}\left\{\frac12
    \int_\Omega\big(|\nabla u_i(x)|^2+|u_i(x)|^2\big)dx-
    \int_\Omega F_i(u_i(x))dx\right\},
\end{equation}
given by the sum of the internal energies of the $k$ densities $u_i$,
each having internal potential $F_i(x,s)=\int_0^s f_i(x,u)du$,
is naturally associated to the system.

Our analysis will highlight how the coexistence of all the densities
is connected to the following minimization problem: \emph{finding
  local minimizers of $J_\Omega(U)$ in the class of of segregated
  states}
\[{\mathcal U}=\left\{U=(u_1,u_2,\dots,u_k) \in
  \big(H^1(\Omega)\big)^k:\, u_i\geq 0,\ u_i\cdot u_j=0\text{ if
  }i\neq j,\text{ a.e. in }\Omega\right\}.\] The problem of the
existence of the \emph{global minimum} of $J_\Omega(U)$ in $\mathcal
U$ was investigated in \cite{ctv-var} under the non-homogeneous
conditions \eqref{dirnonomo}.  As we shall see in Theorem
\ref{p:trivialglobmin}, the global minimizer under homogeneous
boundary conditions is in general trivial, namely a $k$-tuple with all
but one component identically null.  Hence, the only possibility for
finding a \emph{stable} coexistence solution where all the $k$
densities survive, consists in looking for \emph{local} minimizers of
$J_\Omega$, see problem ${\boldsymbol{(P_\e)}}$ below.

Exploiting the variational character of the interaction term in
\eqref{modelLK} and developing a suitable domain perturbation
technique, in this paper we give positive answer to both
questions of minimization of $J_\Omega$ and occurrence of coexistence
states for the system.  Our main result can be summarized as follows:
under suitable assumptions on $f_i$'s ensuring the existence of a
non-degenerate solution to the system on the unperturbed domain
$\Omega_0$ (see \eqref{ND} below), \emph{for all $\Omega_\eps$ close
  enough to $\Omega_0$ and large parameter $\vk$, there exists
  $(u_1^\vk,\dots,u_k^\vk)$ solution to \eqref{modelLK} in
  $\Omega_\eps$, whose limit configuration as $\vk\to\infty$ is a
  segregated coexistence state $(u_1,\dots,u_k)$ with $k$ positive
  components (i.e. each component $u_i\geq 0$ and $u_i$ is strictly
  positive on a set of positive measure), characterized as a local
  minimizer of the free energy~$J_{\Omega_\e}$.}

Before stating rigorously our assumptions and main results, a further
remark is in order.  As observed in \cite{ctv-var} (see also Theorem
\ref{t:extremality}) any $(u_1,\dots,u_k)$ which is a local minimizer
of the free energy $J_\Omega$ on $\mathcal U$, is also a solution of
the following system of distributional inequalities:
\begin{equation}\label{2kvar}
    \begin{cases}
      \displaystyle\int_{\Omega} \big(\nabla u
      _i(x)\nabla\phi(x)+u_i(x) \phi(x)-
      f_i(u _i(x))\phi(x)\big)dx\leq 0,\\[10pt]
      \displaystyle\int_{\Omega}\big(\nabla \widehat u
      _i(x)\nabla\phi(x)+\widehat u_i(x) \phi(x)-\widehat f(\widehat
      u _i(x))\phi(x))dx\geq 0,
    \end{cases}
\end{equation}
$i=1,\dots,k$, for any non-negative $\phi\in H^1(\Omega)$, where we
have denoted $\widehat u_i = u_i-\sum_{h\neq i}u_h$ and $\widehat
f(\widehat u_i) =f(u_i)-\sum_{j\neq i} f_j(u_j)$.  The link between
systems of this form and population dynamics has been pointed out in
\cite{cf,ctv-var,ctv-asymp}: as a matter of fact \emph{all} the
limiting configurations as $\vk\to\infty$ of the solutions to
\eqref{modelLK} are solutions to \eqref{2kvar}. In other words, the possibility of coexistence
of many species ruled out by strong competition is governed by the
system of distributional inequalities \eqref{2kvar}: its independent
study is thus crucial in population dynamics.  In this perspective our
main result can be reformulated in the following way: \emph{ the system of differential
  inequalities \eqref{2kvar} has a solution $(u_1,\dots,u_k)\in
  \mathcal U$ with $k$ positive components.}

\section{Assumptions and main results}

\medskip\noindent {\bf Description of the domain.} We shall work in a
class of smooth non-convex domains $\Omega_\eps$ which generalizes the
dumbbell form with many components as in \cite{d}. 
Let $N\geq 2$ and for
$k\in\N$, let
$$\Omega_0=\Omega^1\cup\Omega^2\cup\dots\cup\Omega^k,$$
where  $\Omega^i\subset \R^N$
are open bounded smooth domains with mutually disjoint closures, i.e.
\begin{equation}\label{unity}
\overline{\Omega^i}\cap\overline{\Omega^j}=\emptyset\qquad\text{ if }i\neq j.
\end{equation}
For any $\e>0$, let
$R_\e\subset(\R^N\setminus \Omega_0)$ be a bounded  measurable set  satisfying
the properties:
\begin{itemize}
\item[(i)]$ |R_\e|\to 0$ as $\e\to 0$
\item[(ii)] $\Omega_0\cup R_\e$ is open and connected
\item[(iii)]  $\partial(\Omega_0\cup R_\e)$ is smooth.
\end{itemize}
Here we denote with $|B|$  the Lebesgue measure of any set $B\subset\R^N$.
Finally  we set
$$\Omega_\eps:=\Omega_0\cup R_\e.$$
%\begin{figure}[h]
%  \centering \subfigure [example of $\Omega_{0}$ with $k=3$]{
%    \begin{psfrags}
%      \psfrag{B1}{$\Omega^1$} \psfrag{B2}{$\Omega^2$}
%      \psfrag{B3}{$\Omega^3$}
%      \includegraphics[width=5.5cm]{triplo1.eps}
%    \end{psfrags}
%  } \quad\quad\quad \subfigure[the set $\Omega_\e$ ] {
%    \includegraphics[width=5.5cm]{triplo2.eps}
%  }
%  \caption{}\label{fig:dd}
%\end{figure}

\medskip\noindent {\bf Assumptions on the nonlinearity.}
For every $i=1,\dots,k$, let $F_i\in C^2(\R)$ with $f_i=F'_i$ satisfying
\begin{itemize}
\item[(F1)] $F_i(0)=0$ and $f_i(0)=0$;
\item[(F2)] there exists $A_i>0$ such that
$\mu_i:=F_i(A_i)-\frac{A_i^2}{2}=
\max_{t\in[0,+\infty)}\big(F_i(t)-\frac{t^2}{2}\big)$;

\smallskip
\item[(F3)] $f'_i(A_i)<1$.
\end{itemize}
Noticeable examples of nonlinearities satisfying (F1)--(F3) are
logistic type functions of the form
$f_i(u)=\lambda u-|u|^{p-1}u$ with $p>1$ and $\lambda>1$.

Assumption (F2) implies that $A_i =f_i(A_i)$ and hence
the constant function $u\equiv A_i$ is a solution to problem
\begin{equation*}
\begin{cases}
-\Delta u+u= f_i(u),&\text{ in } \Omega,\\[5pt]
\;\;\;\dfrac{\partial u}{\partial \nu}=0,&\text{ on } \partial\Omega,
\end{cases}
\end{equation*}
in any open smooth domain $\Omega$; moreover  $u\equiv A_i$ minimizes
the internal energy
$$
\int_\Omega\left(\frac12
|\nabla u_i(x)|^2+\frac12|u_i(x)|^2- F_i(u_i(x))\right)dx.
$$

\medskip\noindent
Let us denote $w_i=A_i\alchi_{\Omega^i}$,
$i=1,\dots,k$, $W=(w_1,w_2,\dots,w_k)\in \big(H^1(\Omega_0)\big)^k$,
and set
\begin{equation}\label{e:mu}
\mu=\sum_{i=1}^k\left\{ \frac12\int_{\Omega_0}\big(|\nabla
  w_i|^2+|w_i|^2\big)\,dx-\int_{\Omega_0} F_i(w_i)\,dx\right\}
=-\sum_{i=1}^k\mu_i|\Omega^i|.
\end{equation}
Assumption (F3) implies the following non-degeneracy property:
for all $i$ and $u\in H^1(\Omega^i)$
\begin{equation}\label{ND}
  \int_{\Omega^i} (|\nabla u|^2+|u|^2-f_i'(w_i)u^2)dx
  \geq\nu\int_{\Omega^i}(|\nabla u|^2+|u|^2) dx,
\end{equation}
where $\nu:=\min_{i=1,\dots,k}\{1,1-f_i'(A_i)\}>0$.

\medskip
We are now going to describe the  main results of the present paper,
 starting from the following optimal partition problem.

\medskip\noindent
{\bf Problem ${\boldsymbol{(P_\e)}}$}.
Find \emph{nontrivial} local minimizers of the functional
\begin{align*}
&J_{\Omega_\e}:\big(H^1(\Omega_\e)\big)^k\to(-\infty,+\infty],\\
&  J_{\Omega_\e}(U)=\sum_{i=1,\dots, k}\left\{\frac12
    \int_{\Omega_\e}\big(|\nabla u_i(x)|^2+|u_i(x)|^2\big)dx-
    \int_{\Omega_\e} F_i(u_i(x))dx\right\},
\end{align*}
among $k$-tuples $U=(u_1,u_2,\dots,u_k)$ belonging to the class
\[
{\mathcal U}_\e=\left\{U=(u_1,u_2,\dots,u_k) \in
  \big(H^1(\Omega_\e)\big)^k:\, u_i\geq 0,\ u_i\cdot u_j=0\text{ if }i\neq
  j,\text{ a.e. in }\Omega_\e\right\}.
\]
By \emph{nontrivial} we mean that no component $u_i$
of the solution $U$ can be null, i.e. $u_i\not\equiv 0$ for all
$i=1,\dots,k$.
As stated in the introduction, we shall prove that in any connected domain and
for a wide class of $F_i$'s
including logistic-type nonlinearities, any
\emph{global minimizer} of the free energy \eqref{freeenergy} is indeed trivial.

\begin{Proposition}\label{p:trivialglobmin}
  Let $\Omega\subset\R^N$ be a connected open domain and $F_i\in
  C^2(\R)$ satisfy  $F_i(0)=0$ and {\em{(F2)}}.
Then the infimum
\[
\lambda:=\inf_{U\in \mathcal U}J_\Omega(U)
\]
is achieved by $U_0=(u_1^0,\dots,u_k^0)$
with $u_{i_0}^0\equiv A_{i_0}$ and $u_i^0\equiv 0$ for $i\neq i_0$, and
\[\lambda=
-\mu_{i_0}|\Omega|,
\]
where $\mu_{i_0}=\max_{i\in \{1,\dots,k\}}\mu_i$.  Furthermore, any
$k$-tuple achieving $\lambda$ has all but one component identically
null.
\end{Proposition}
In view of the above proposition, there is no hope to find
nontrivial solutions to ${\boldsymbol{(P_\e)}}$ by global
minimization.  On the contrary, by studying $J_{\Omega_\e}$ \emph{near $W$},
we can find positive answer to the problem. To this aim let us
denote by
\[
B^\delta_\e(W):=\left\{U\in \big(H^1(\Omega_\e)\big)^k:
  \,\|U-W\|_{(H^1(\Omega_0))^k}<\delta\right\}
\]
the set of $k$-tuples $U$ whose restriction to $\Omega_0$ is close
within $\delta>0$ to $W$, with respect to the $H^1$ norm
$\|V\|^2_{(H^1(\Omega_0))^k}= \sum_{i=1}^k\|v_i\|^2_{H^1(\Omega_0)}$.
Notice that, if $U=(u_1,\dots,u_k)\in B^\delta_\e(W)$, then each  $u_i$
satisfies  $\int_{\Omega_i}|u_i-A_i|^2\leq \delta^2$.
Hence, if
\begin{equation}\label{eq:nontrivi}
\delta^2<A_i^2 |\Omega_i|,\qquad i=1,\dots,k
\end{equation}
then $u_i\not\equiv0$. Henceforward, $\delta$ will be supposed
to satisfy \eqref{eq:nontrivi}, thus ensuring that any $U\in
B^\delta_\e(W)$ is nontrivial.
\begin{Theorem}\label{t:partizione}
Assume that (F1)--(F3) hold and let
\[
\lambda_{\e}^{\delta}:=\inf_{U\in{\mathcal
      U}_\e\cap {B^\delta_\e(W)}} J_{\Omega_\e}(U).
\]
Then, there exists $\delta>0$ such that, for every $\e$ sufficiently
small, $\lambda_{\e}^{\delta}$ is achieved by a $k$-tuple
$U_{\e}=(u_1^{\e},\dots,u_k^{\e})$ with $0\leq u_i^\e\leq A_i$ a.e.
in $\Omega_\e$ and $u_i^\e\not\equiv 0$ for all $i=1,\dots,k$.
\end{Theorem}
The proof of Theorem \ref{t:partizione} will be obtained through a
careful analysis of the solutions to the original competitive system
\eqref{modelLK}, as the parameter $\vk$ of the interspecific
competition grows. Our main result reads as follows:
\begin{Theorem}\label{t:sistemi}
  Assume that (F1)--(F3) hold. Then, there exists $\delta>0$ such
  that, for every $\e$ sufficiently small and and $\vk>0$ sufficiently
  large, system \eqref{modelLK} coupled with \eqref{neuomo} in
  $\Omega_\e$ admits a solution
  $U^{\e,\vk}=(u_1^{\e,\vk},\dots,u_k^{\e,\vk})\in B^\delta_\e(W)$
  with the following properties:
  \begin{enumerate}
  \item $u_i^{\e,\vk}\not\equiv 0$ for all $i=1,\dots,k$.
  \item $0\leq u_i^{\e,\vk}\leq A_i$ a.e. in $\Omega_\e$ for all
    $i=1,\dots,k$.
  \item There exists $V^\e=(v_1^{\e},...,v_k^{\e})\in \mathcal
    U_\e\cap B^\delta_\e(W)$ such that $v_i^\e\not\equiv 0$ for every
    $i$ and, up to subsequences, $U^{\e,\vk}\to V^\e$ strongly in
    $(H^1(\Omega_\e))^k$ as $\vk\to+\infty$.  Furthermore, $V^\e$ is a
    local minimizer of $J_{\Omega_\e}$, namely
    $J_{\Omega_\e}(V^\e)=\lambda_\e^\delta$.
\end{enumerate}
\end{Theorem}
\noindent The proof of our results relies on the minimization on
the whole $B^\delta_\e(W)$ of an auxiliary functional obtained by
penalizing the internal energy $J_{\Omega_\e}$ with a positive competition
term. More precisely we shall consider a suitable modification of the
functional
\begin{multline*}
  \sum_{i=1}^k\left\{
    \frac12\int_{\Omega_\e}\big(|\nabla
    u_i(x)|^2+|u_i(x)|^2\big)\,dx-\int_{\Omega_\e}
   F_i(u_i(x))\,dx\right\}
  +\k \sum_{\substack{i,j=1\\i\neq j}}^k \int_{\Omega_\e}
  u_i(x)^2u_j(x)^2\,dx
\end{multline*}
defined on $(H^1(\Omega_\e))^k$, see $I_{\eps,\vk}$ in \eqref{eq:Iek}
below. Due to the variational character of the competition term in
\eqref{modelLK}, by standard Critical Point Theory, any local
minimizer of the above function is a (weak) solution to the original
system. Section \ref{s:Iekappa} is devoted to the search for a local
minimizer of $I_{\eps,\vk}$ in $B^\delta_\e(W)$ and requires the main
technical effort of the paper. By developing a domain perturbation
argument based on the nondegeneracy condition \eqref{ND}, we shall
succeed in proving the existence of a minimizer in small perturbations
of the domain $\Omega_0$, for large values of the competition
parameter $\vk$.  In this way, we directly obtain the existence of a
positive solution to the competitive system, at any fixed $\vk$, see
Section~\ref{sec:competitive-systems}. In the subsequent Section
\ref{sec:asymptotic-analysis} we perform the asymptotic analysis of
these solutions as the competition parameter $\vk\to\infty$, showing
that the steady states segregate in a nontrivial limit configuration
$V^\e$. The comparison between the minimal energy levels of
$I_{\eps,\vk}$ and $J_{\Omega_e}$ will allow proving that $V^\e$
indeed solves problem ${\boldsymbol{(P _\e)}}$ on $B^\delta_\e(W)$.
This concludes the proof of Theorem \ref{t:sistemi} and, in turn, that
of Theorem \ref{t:partizione}. In the last part of Section
\ref{sec:asymptotic-analysis}, we show that any solution to the
optimal partition problem ${\boldsymbol{(P _\e)}}$ satisfies some
extremality conditions in the form of differential inequalities
\eqref{2kvar}.  Finally, in the last section we derive some
consequences of this fact, and outline further developments of the
subject.

\section{A variational problem}\label{s:Iekappa}
Aim of this section is to study the minimization of a suitable
functional on $(H^1(\Omega_\e))^k$, which will reveal to be strongly
related both to problem ${\boldsymbol{(P _\e)}}$ and to the original
competitive system. The functional is defined as follows:
\begin{multline}\label{eq:Iek}
  I_{\e,\k}(U)=\sum_{i=1}^k\left\{
    \frac12\int_{\Omega_\e}\big(|\nabla
    u_i(x)|^2+|u_i(x)|^2\big)\,dx-\int_{\Omega_\e}
    \widetilde F_i(u_i(x))\,dx\right\}\\
  +\k \sum_{\substack{i,j=1\\i\neq j}}^k \int_{\Omega_\e}
  G_i(u_i(x))G_j(u_j(x))\,dx
\end{multline}
where
\[
\widetilde F_i(t)=\begin{cases}
0,\quad &\text{if }t\leq 0,\\
F_i(t),\quad &\text{if }0\leq t\leq  A_i,\\
A_it+F_i(A_i)-A_i^2,\quad &\text{if }t\geq A_i,
\end{cases}
\]
and
\[
G_i(t)=\begin{cases}
t^2,\quad &\text{if }|t|\leq A_i,\\
2A_i|t|-A_i^2,\quad &\text{if }|t|> A_i.
\end{cases}
\]
Notice that  $I_{\e,\k}\in C^1\big((H^1(\Omega_\e))^k,\R\big)$.
Aim of this section is to prove
\begin{Theorem}\label{t:auxiliary}
Assume that (F1)--(F3) hold and let
\[
c_{\e,\vk}:=\inf_{U\in{\mathcal
      U}_\e\cap {B^\delta_\e(W)}} I_{\e,\vk}(U).
\]
Then, there exists $\delta>0$ such that, for every for $\e>0$
sufficiently small and $\vk>0$ sufficiently large, $c_{\e,\vk}$ is
achieved by a $k$-tuple $U^{\e,\vk}=(u_1^{\e,\vk},\dots,u_k^{\e,\vk})$
with $0\leq u_i^{\e,\vk}\leq A_i$ a.e.  in $\Omega_\e$ and
$u_i^{\e,\vk}\not\equiv 0$ for all $i=1,\dots,k$.
\end{Theorem}

The first step in this direction consists in proving that the minimum
is achieved on the closure of $B^\delta_\e(W)$, namely the set
\[
\overline{B^\delta_\e(W)}:=\left\{U\in \big(H^1(\Omega_\e)\big)^k:
  \,\|U-W\|_{(H^1(\Omega_0))^k}\leq \delta\right\}.
\]
\begin{Lemma}\label{l:Lambdaeps}
  For every $\delta$ satisfying \eqref{eq:nontrivi}, $\e
  \in(0,1)$, and $\vk>0$, the infimum
  \[
  \Lambda_{\e,\vk}=\inf_{U\in \overline{B^\delta_\e(W)}} I_{\e,\vk}(U)
  \]
  is achieved by a $k$-tuple
  $U^{\e,\vk}=(u_1^{\e,\vk},\dots,u_k^{\e,\vk})$ where
  $u_i^{\e,\vk}\not\equiv 0$ and
\begin{equation}\label{eq:apriori1sis}
  0\leq u_i^{\e,\vk}(x)\leq A_i\quad \text{for a.e. }x\in\Omega_\e.
\end{equation}
\end{Lemma}
\begin{pf}
  We first observe that $\frac12 t^2-\widetilde F_i(t)\geq \frac12
  A_i^2-F_i(A_i)$ for all $t\in\R$, hence, being the coupling term
  nonnegative, for all $U=(u_1,\dots,u_k)\in
  \big(H^1(\Omega_\e)\big)^k$
\[
I_{\e,\vk}(U)\geq \sum_{i=1}^k\int_{\Omega_\e}\bigg(\frac12
|u_i|^2-\widetilde F_i(u_i)\bigg)dx\geq
\sum_{i=1}^k\bigg(\frac{A_{i}^2}2-F_{i}(A_{i})\bigg)|\Omega_\e|,
\]
and hence $\Lambda_{\e,\vk}>-\infty$.  Let
$\big\{U_n=(u_1^n,\dots,u_k^n)\big\}_{n\in\N}$ be a minimizing
sequence, i.e. $ U_n\in{\overline{B^\delta_\e(W)}}$ and
$\lim_{n\to+\infty}I_{\e,\vk}(U_n)=\Lambda_{\e,\vk}$.  We notice
that, by definition of $\widetilde F_i$ and the fact that $w_i\geq
0$ a.e., we can choose $U_n$ such that $u_i^n\geq 0$ a.e. in
$\Omega_\e$ for all $i=1,\dots,k$ (otherwise we take
$((u_1^n)^+,\dots,(u_k^n)^+)$ with $(u_i^n)^+:=\max\{u_i^n,0\}$ as
a new minimizing sequence). Letting $V_n=(v_1^n,\dots,v_k^n)$ with
$v_i^n=\min\{u_i^n,A_i\}$, it is easy to verify that $V_n\in
\overline{B^\delta_\e(W)}$ and $I_{\e,\vk}(V_n)\leq
I_{\e,\vk}(U_n)$.  Then also $\big\{V_n\big\}_{n\in\N}$ is a
minimizing sequence.

Since $\big\{V_n\big\}_{n\in\N}$ is a minimizing sequence and it is
uniformly bounded , it is easy to realize that
$\big\{V_n\big\}_{n\in\N}$ is bounded in $\big(H^1(\Omega_\e)\big)^k$,
hence there exists a subsequence, still denoted as
$\big\{V_n\big\}_{n\in\N}$, which converges to some
$V=(v_1,\dots,v_k)\in \big(H^1(\Omega_\e)\big)^k$ weakly in
$\big(H^1(\Omega_\e)\big)^k$, strongly in $\big(L^2(\Omega_\e)\big)^k$
and a.e. in $\Omega_\e$.  A.e.  convergence implies that $0\leq
v_i\leq A_i$ a.e. in $\Omega_\e$, while weakly lower semi-continuity
implies that $V\in \overline{B^\delta_\e(W)}$.  From $0\leq v_i\leq
A_i$ and the Dominated Convergence Theorem, it follows that
\begin{align*}
  &\lim_{n\to+\infty}\int_{\Omega_\e}\widetilde
  F_i(v_i^n(x))\,dx=\int_{\Omega_\e}\widetilde F_i(v_i(x))\,dx,\\
  &\lim_{n\to+\infty}\int_{\Omega_\e} G_i(v_i^n(x))G_j(v_j^n(x))\,dx
  =\int_{\Omega_\e} G_i(v_i(x))G_j(v_j(x))\,dx,
\end{align*}
for every $i,j=1,\dots,k$, which, together with  lower
semi-continuity, yields
\[
\Lambda_{\e,\vk}\leq I_{\e,\vk}(V)\leq
\liminf_{n\to+\infty}I_{\e,\vk}(V_n)=
\lim_{n\to+\infty}I_{\e,\vk}(V_n)=\Lambda_{\e,\vk},
\]
thus proving that $V$ attains $\Lambda_{\e,\vk}$.

 Finally, if
$v_i\equiv 0$ in $\Omega_i$ then
$\|v_i-A_i\|_{H^1(\Omega^i)}^2=\int_{\Omega^i} A_i^2 dx\leq
\delta^2$, in contradiction with the choice of $\delta$ as in \eqref{eq:nontrivi}.
\end{pf}

A major effort is now needed to show that the minimum provided by
Lemma \ref{l:Lambdaeps} indeed belongs to the open set
$B_\delta^\e(W)$.  The crucial ingredient in this direction
consists in providing suitable estimates of the minimal level
$\Lambda_{\e,\vk}$, which require the following technical lemma.
\begin{Lemma}\label{l:continuity}
  For every $\eta>0$ there exists $\delta_{\eta}>0$ such that if
  $U\!=\!(u_1,\dots,u_k)\in \overline{B^{\delta_\eta}_\e\!(W)}$ and $|u_i(x)|\leq
  A_i$ for a.e. $x\in\Omega_0$ and for all $i=1,\dots,k$, then
  \begin{equation}\label{eq:continuity}
  \sum_{i=1}^k\int_{\Omega_0}\bigg[F_i(u_i)-F_i(w_i)-
  f_i(w_i)(u_i-w_i)-\frac12 f_i'(w_i)(u_i-w_i)^2\bigg]\,dx\leq\eta
  \|U-W\|^2_{(H^1(\Omega_0))^k}.
\end{equation}
\end{Lemma}
\begin{pf}
We have
\begin{align*}
  &\int_{\Omega_0}\bigg[F_i(u_i)-F_i(w_i)- F_i'(w_i)(u_i-w_i)-
  \frac12 F_i''(w_i)(u_i-w_i)^2\bigg]\,dx\\
  &=\int_{\Omega_0}\!\bigg[\!
  \int_0^1\!\bigg[\!\bigg(\!\frac{d}{dt}F_i(t\,u_i+(1-t)w_i)\!\bigg)\!-
  F_i'(w_i)(u_i-w_i)- t\, F_i''(w_i)(u_i-w_i)^2\bigg]dt\bigg]dx\\
&=\int_{\Omega_0}\!\bigg[\!
  \int_0^1\!\big[F_i'(t\,u_i+(1-t)w_i)-
  F_i'(w_i)- t\, F_i''(w_i)(u_i-w_i)\big](u_i-w_i)\,dt\bigg]dx
\\
&=\int_{\Omega_0}\!\bigg[\!
  \int_0^1\!\!
  \bigg(\int_0^1\!\!\bigg(\!\frac{d}{ds}
F_i'\big(s(t\,u_i+(1-t)w_i)+(1-s)w_i\big)\bigg)ds\\
&\hskip7cm-
  t\, F_i''(w_i)(u_i-w_i)\!\bigg)(u_i-w_i)\,dt\bigg]dx.
\end{align*}
Hence, by H\"older's inequality,
\begin{align*}
  \int_{\Omega_0}&\Big[F_i(u_i)-F_i(w_i)-
  F_i'(w_i)(u_i-w_i)-
  \frac12 F_i''(w_i)(u_i-w_i)^2\Big]\,dx\\
  &\leq \int_{\Omega_0}\!\bigg[\!  \int_0^1\!  \bigg(\int_0^1\!\big(
  F_i''\big(st(u_i-w_i)+w_i\big)-
  F_i''(w_i)\big)t(u_i-w_i)^2ds\bigg)\,dt\bigg]\,dx\\
  &\leq
  \|u_i-w_i\|_{L^{p}(\Omega_0)}^2\iint_{(0,1)\times(0,1)}t\|F_i''\big(st(u_i-w_i)
  +w_i\big)-F_i''(w_i)\|_{L^{\frac{p}{p-2}}(\Omega_0)}\,ds\,dt,
\end{align*}
where $p=2^*$ for $N\geq3$ and $p\in(2,+\infty)$ for $N=2$.  The
conclusion follows now from Sobolev's embeddings and the
continuity of the operator
\begin{align*}
   F_i'':\big\{v\in H^1(\Omega_0):\,|v(x)|\leq 3A_i\big\}&\to
  L^{\frac{p}{p-2}}(\Omega_0),\\
   v&\mapsto F_i''(v),
\end{align*}
which can be easily proved using  the Dominated Convergence
Theorem.
\end{pf}

\begin{remark}\label{rem:ipotesi}
 According to Lemma \ref{l:continuity},
 besides (\ref{eq:nontrivi}) from now on we assume
\begin{equation*}
0<\delta\leq \delta_0
\end{equation*}
with $\delta_0$ small enough in such a way that inequality
(\ref{eq:continuity}) with $\eta=\min\{\frac{\nu}{4},\frac18\}$
holds for all functions $U=(u_1,\dots,u_k)\in
\overline{B^{\delta_0}_\e(W)}$ satisfying $|u_i(x)|\leq A_i$ a.e.
in $\Omega_0$. We also require that $\delta_0\leq A_i^2/4$ and finally
that condition \eqref{eq:sobo} in Lemma \ref{l:stima2} is
satisfied.
\end{remark}

\noindent By exploiting the separation of the $\Omega^i$'s as in \eqref{unity}, for every $i=1,\dots,k$, 
we can construct test functions $\varphi^i\in H^1(\R^N)$
satisfying
\begin{equation}\label{eq:apriori2}
0\leq\varphi^i(x)\leq A_i \quad \text{a.e. in }\R^N,
\end{equation}
$\varphi_i(x)=0$ for all $x\in \Omega_0\setminus \Omega_i$,
$\varphi_i(x)=A_i$ if $x\in\Omega_i$, and $\varphi_i\cdot\varphi_j=0$
a.e. in $\R^N$ if $i\neq j$.  
This allows us to provide an estimate from above of the value $\Lambda_{\e,\vk}$
in terms of the total free-energy of $W$.
\begin{Lemma}\label{l:stima1}
  For every $\e\in(0,1)$, there exists $\tau_\e$ such that $\tau_\e\to
  0$ as $\e\to 0$ and, for all $\vk>0$,
  \[
  \Lambda_{\e,\vk}\leq \mu+\tau_\e,
  \]
with $\mu$ given by \eqref{e:mu}.
\end{Lemma}
\begin{pf}
Let $\varphi_\e^i\in
H^1(\Omega_\e)$ be the restriction of $\varphi_i$ to $\Omega_\e$.
Notice that $(\varphi_\e^1,\varphi_\e^2,\dots,\varphi_\e^k)\in
  \overline{B^\delta_\e(W)}$ and that $\varphi_\e^i\cdot\varphi_\e^j\equiv
  0$ if $i\neq j$. Hence we have 
\begin{align*}
  \Lambda_{\e,\vk}&\leq  I_{\e,\vk}(\varphi_\e^1,\varphi_\e^2,\dots,\varphi_\e^k)\\
  &=\mu+ \sum_{i=1}^k\left\{ \frac12\int_{R_\e}\big(|\nabla
    \varphi_\e^i(x)|^2+|\varphi_\e^i(x)|^2\big)\,dx-\int_{R_\e}
    F_i(\varphi_\e^i)\,dx\right\}\\
    &= \mu+\tau_\e,
\end{align*}
where  \begin{equation*}
\tau_\e=\sum_{i=1}^k\left\{ \frac12\int_{R_\e}\big(|\nabla
    \varphi^i(x)|^2+|\varphi^i(x)|^2\big)\,dx-\int_{R_\e}
    F_i(\varphi^i)\,dx\right\}.
\end{equation*}
Since $|R_\e|\to 0$ as $\e\to 0$, then $\tau_\e\to 0$,
proving the stated estimate.
\end{pf}
\begin{Lemma}\label{l:stima2}
  For every $\e\in(0,1)$, there exists $\sigma_\e$ such that
  $\sigma_\e\to 0$ as $\e\to0$ and
  \[
  \|U^{\e,\vk}-W\|^2_{(H^1(\Omega_0))^k}\leq \sigma_\e
  \]
for every $\vk>\max\limits_{i\neq j }\frac{2f_i'(0 )}{A_j^2}$.
\end{Lemma}
\begin{pf}
  From (\ref{eq:apriori1sis}), we can write
  $\Lambda_{\e,\vk}=I^1_{\e,\vk}+I^2_{\e,\vk}$ where
\begin{align*}
  &I^1_{\e,\vk}=\sum_{i=1}^k\bigg\{
  \frac12\int_{\Omega_0}\!\!\big(|\nabla
  u_i^{\e,\vk}|^2+|u_i^{\e,\vk}|^2\big)\,dx-\!\!\int_{\Omega_0}\!\!
  F_i(u_i^{\e,\vk})\,dx+ \k \sum_{j\neq i} \int_{\Omega_0}\!\!
  (u_i^{\e,\vk})^2(u_j^{\e,\vk})^2\,dx\bigg\}\\
  &I^2_{\e,\vk}=\sum_{i=1}^k\bigg\{ \frac12\int_{R_\e}\!\!\big(|\nabla
  u_i^{\e,\vk}|^2+|u_i^{\e,\vk}|^2\big)\,dx-\!\!\int_{R_\e}\!\!
  F_i(u_i^{\e,\vk})\,dx + \k \sum_{j\neq i} \int_{R_\e}\!\!
(u_i^{\e,\vk})^2(u_j^{\e,\vk})^2\,dx\bigg\}.
\end{align*}
Since by assumption $-\Delta w_i+w_i=f_i(w_i)$ in $\Omega_0$, we
can write each term in $I^1_{\e,\vk}$ as follows
\begin{align*}
  \frac12\int_{\Omega_0}&\big(|\nabla u_i^{\e,\vk} |^2+|u_i^{\e,\vk}
  |^2\big)\,dx-\int_{\Omega_0} F_i(u_i^{\e,\vk} )\,dx+ \k \sum_{j\neq
    i} \int_{\Omega_0}
  (u_i^{\e,\vk})^2(u_j^{\e,\vk})^2\,dx\\
  &=\frac12\int_{\Omega_0}\big(|\nabla w_i |^2+|w_i
  |^2\big)\,dx-\int_{\Omega_0}
  F_i(w_i )\,dx\\
  &\quad+ \frac12\int_{\Omega_0}\big(|\nabla (u_i^{\e,\vk}-w_i)
  |^2+|(u_i^{\e,\vk}-w_i) |^2\big)\,dx
  -\int_{\Omega_0} \big( F_i(u_i^{\e,\vk} )-F_i(w_i )\big)\,dx\\
  &\quad+ \int_{\Omega_0}\big(\nabla w_i \cdot\nabla(u_i^{\e,\vk}-w_i)
  +w_i (u_i^{\e,\vk}-w_i) \big)\,dx+ \k \sum_{j\neq i} \int_{\Omega_0}
  (u_i^{\e,\vk})^2(u_j^{\e,\vk})^2\,dx\\
  &=-\mu_i|\Omega^i| +\alpha^1_{\e,\vk,i}+\alpha^2_{\e,\vk,i}\\
  &\quad -\int_{\Omega_0} \big( F_i(u_i^{\e,\vk} )-F_i(w_i )- f_i(w_i
  )(u_i^{\e,\vk}-w_i)-\frac12f_i'(w_i )(u_i^{\e,\vk}-w_i)^2 \big)\,dx.
\end{align*}
where
\[
\alpha^1_{\e,\vk,i}=\frac12\|u_i^{\e,\vk}-w_i\|_{H^1(\Omega^i)}^2-
  \frac12 \int_{\Omega^i}f_i'(A_i )(u_i^{\e,\vk}-w_i)^2 + \k
  \sum_{j\neq i} \int_{\Omega^i}
  (u_i^{\e,\vk})^2(u_j^{\e,\vk})^2\,dx
\]
and
\[
\alpha^2_{\e,\vk,i}=\frac12 \sum_{j\neq i}\int_{\Omega^j}\bigg(
|\nabla u_i^{\e,\vk}|^2+|u_i^{\e,\vk}|^2- \Big[f_i'(0 )- 2\k
  \sum_{h\neq i}(u_h^{\e,\vk})^2\Big]|u_i^{\e,\vk}|^2\bigg)\,dx.
\]
From (\ref{ND}) it follows that
\begin{equation}\label{eq:6}
\alpha^1_{\e,\vk,i}\geq
\frac\nu2\|u_i^{\e,\vk}-w_i\|_{H^1(\Omega^i)}^2.
\end{equation}
On the other hand, from H\"older's and Sobolev's inequalities it
follows that
\[
\alpha^2_{\e,\vk,i}\geq\frac12 \sum_{j\neq i}\int_{\Omega^j}\bigg(
\Big(1-\Big\|\Big[f_i'(0 )- 2\k
  \sum_{h\neq i}(u_h^{\e,\vk})^2\Big]^+\Big\|_{L^{\frac{p}{p-2}}(\Omega^j)}
S_{p,j}^{-1}\Big) |\nabla
u_i^{\e,\vk}|^2+|u_i^{\e,\vk}|^2\bigg)\,dx.
\]
where $p=2^*$ for $N\geq3$ and $p\in(2,+\infty)$ for $N=2$, and
$S_{p,j}$ is the best constant in the Sobolev embedding
$H^1(\Omega^j)\hookrightarrow L^p(\Omega^j)$. Let us denote
\[
A^\delta_{\vk,j}=\{x\in\Omega^j:|u_j^{\e,\vk}-A_j|^2>\delta\}.
\]
Hence
\[
\delta^2\geq\int_{\Omega^j}|u_j^{\e,\vk}-A_j|^2\,dx\geq \delta
|A^\delta_{\vk,j}|
\]
and then $|A^\delta_{\vk,j}|<\delta$. In particular, if $\delta$
is such that
\begin{equation}\label{eq:sobo}
\delta^{\frac{p-2}{p}}|(f_i'(0 ))^+|<\frac{S_{p,j}}2,
\end{equation}
there holds
\begin{equation}\label{eq:4}
\Big\|\Big[f_i'(0 )- 2\k \sum_{h\neq
  i}(u_h^{\e,\vk})^2\Big]^+\Big\|_{L^{\frac{p}{p-2}}(A^\delta_{\vk,j})}
<\frac{S_{p,j}}2.
\end{equation}
In $\Omega^j\setminus A^\delta_{\vk,j}$, there holds
$u_j^{\e,\vk}>A_j-\sqrt \delta>\frac{A_j}{2}$ for $\delta$ small
as in Remark \ref{rem:ipotesi}. Then, if $\vk>{2f_i'(0
)}/{A_j^2}$,
\begin{equation}\label{eq:5}
f_i'(0 )- 2\k \sum_{h\neq i}(u_h^{\e,\vk})^2<0 \quad \text{in
}\Omega^j\setminus A^\delta_{\vk,j}.
\end{equation}
Collecting (\ref{eq:4}) and (\ref{eq:5}), we deduce that, for $\vk
>\frac{2f_i'(0 )}{A_j^2}$,
\[
\Big\|\Big[f_i'(0 )- 2\k
  \sum_{h\neq i}(u_h^{\e,\vk})^2\Big]^+\Big\|_{L^{\frac{p}{p-2}}(\Omega^j)}
S_{p,j}^{-1}<\frac12,
\]
and therefore
\begin{equation}\label{eq:7}
\alpha^2_{\e,\vk,i}\geq \frac14 \sum_{j\neq
i}\|u_i^{\e,\vk}-w_i\|_{H^1(\Omega^j)}^2.
\end{equation}
From (\ref{eq:6}) and (\ref{eq:7}), we obtain that
\[
\alpha^1_{\e,\vk,i}+\alpha^2_{\e,\vk,i}\geq
\min\Big\{\frac\nu2,\frac14\Big\}
\|u_i^{\e,\vk}-w_i\|_{H^1(\Omega_0)}^2.
\]
By Lemma \ref{l:continuity} and Remark \ref{rem:ipotesi} we have
that
\begin{multline*}
    \sum_{i=1}^k\int_{\Omega_0}\big[F_i(u_i^{\e,\vk} )-F_i(w_i )- f_i(w_i
  )(u_i^{\e,\vk}-w_i)-\frac12f_i'(w_i )(u_i^{\e,\vk}-w_i)^2 \big]\,dx\\
    \leq
\min\Big\{\frac\nu4,\frac18\Big\}
\|U^{\e,\vk}-W\|^2_{(H^1(\Omega_0))^k}.
  \end{multline*}
Hence
\begin{equation}\label{eq:i1}
  I^1_{\e,\vk}\geq \mu+\min\Big\{\frac\nu4,\frac18\Big\}
\|U_{\e}-W\|^2_{(H^1(\Omega_0))^k}.
\end{equation}
On the other hand, $I_{\e,\vk}^2$ can be promptly estimated by
\begin{align}\label{eq:i2}
I^2_{\e,\vk}\geq -|R_\e|\sum_{i=1}^k\mu_i
\end{align}
with $\mu_i$ as in (F2).
Combining inequalities \eqref{eq:i1} and \eqref{eq:i2}, it follows
that
\begin{equation}\label{eq:sti2}
\Lambda_{\e,\vk}\geq
\mu+\eta\|U_{\e}-W\|^2_{(H^1(\Omega_0))^k}-|R_\e|\sum_{i=1}^k\mu_i,
\end{equation}
where $\eta=\min\{\frac{\nu}{4},\frac18\}>0$. From Lemma
\ref{l:stima1} and (\ref{eq:sti2}), we infer that
$\|U_{\e}-W\|^2_{(H^1(\Omega_0))^k}\leq \sigma_\e$ with
$\sigma_\e=\frac{1}{\eta}(\tau_\e+|R_\e|\sum_{i=1}^k\mu_i)$, concluding the proof.
\end{pf}

\medskip\noindent
\begin{pfn}{Theorem \ref{t:auxiliary}}
  In order to conclude the proof of the theorem, it is sufficient to
  consider $U^{\e,\vk}$ provided by Lemma \ref{l:Lambdaeps}. If
  $\vk$ is large enough, we can apply Lemma \ref{l:stima2} and we
  infer that $U^{\e,\vk}\in B^\delta_\e(W)$ provided $\e$ is
  sufficiently small, and hence $U^{\e,\vk}$ attains
  $c_{\e,\vk}=\Lambda_{\e,\vk}$, i.e. it is a local minimizer of
  $I_{\e,\vk}$ on the open set $B^\delta_\e(W)$ with all the required
  properties.
\end{pfn}

\section{Competitive systems}\label{sec:competitive-systems}
In this section we prove the existence of solutions to the competitive system
\begin{equation}
\label{e:nostro}
\begin{cases}
  -\Delta u_i+u_i=f_i(u_i)- 2\varkappa u_i\sum_{j\neq i}u_j^2,
  &\text{ in }\Omega_\e,\\
  \displaystyle\frac{\partial u_i}{\partial \nu}=0,&\text{ on
  }\partial\Omega_\e,
 \end{cases}
\end{equation}
for $i=1,\dots,k$.
\begin{Theorem}\label{t:sistemi2}
  There exists $\delta>0$ such that for $\e>0$ sufficiently small and
  $\vk>0$ sufficiently large, system \eqref{e:nostro} admits a
  solution $U^{\e,\vk}=(u_1^{\e,\vk},\dots,u_k^{\e,\vk})\in
  B^\delta_\e(W)$ such that, for all $i=1,\dots,k$,
  $u_i^{\e,\vk}\not\equiv 0$ and
\begin{equation}\label{e:apr}
0\leq   u_i^{\e,\vk}\leq A_i\quad\text{a.e. in }\Omega_\e.
\end{equation}
\end{Theorem}
\begin{proof}
By standard Critical
Point Theory, see e.g. \cite{am}, the critical points of
$I_{\e,\k}$ on $(H^1(\Omega_\e))^k$ give rise to weak (and by
regularity classical) solutions to
\begin{equation}
\label{e:truncated}
\begin{cases}
  -\Delta u_i+u_i=\widetilde f_i(u_i)- \varkappa g_i(u_i)\sum_{j\neq i}G_j(u_j),
  &\text{ in }\Omega_\e,\\
  \displaystyle\frac{\partial u_i}{\partial \nu}=0,&\text{ on
  }\partial\Omega_\e,
 \end{cases}
\end{equation}
where
\[
\widetilde f_i(t)=\begin{cases}
0,\quad &\text{if }t\leq 0,\\
f_i(t),\quad &\text{if }0\leq t\leq  A_i,\\
A_i,\quad &\text{if }t\geq A_i,
\end{cases}
\]
and
\[
g_i(t)=\begin{cases}
2t,\quad &\text{if }|t|\leq A_i,\\
2A_i\sgn(t),\quad &\text{if }|t|> A_i.
\end{cases}
\]
Notice that a solution to \eqref{e:truncated} satisfying \eqref{e:apr}
is also a solution of \eqref{e:nostro}.  Now the proof of the theorem
immediately follows by considering $U^{\e,\vk}\in B^\delta_\e(W)$ as
in Theorem \ref{t:auxiliary}; since it is a local minimizer of
$I_{\e,\vk}$, it is a free critical point of $I_{\e,\vk}$ and hence
solves \eqref{e:truncated}. By the validity of \eqref{eq:apriori1sis}
we finally deduce that $U^{\e,\vk}$ is actually a solution to
\eqref{e:nostro}, thus completing the proof.

\end{proof}

\section{The optimal partition problem}\label{sec:asymptotic-analysis}

In this section we deal with problem ${\boldsymbol{(P_\e)}}$,
namely we look for local minimizers of the free energy on
segregated states. The localization of the problem is essentially
motivated by the fact that
 any global minimizer of the free energy in a connected
domain is trivial, as stated in Proposition
\ref{p:trivialglobmin}, the proof of which is given below.

\medskip\noindent
\begin{pfn}{Proposition \ref{p:trivialglobmin}}
By a direct computation, for any $U=(u_1,\dots,u_k) \in{\mathcal
U}$
\begin{align}\label{eq:1}
  J_{\Omega}(U)&\geq \sum_{i=1}^k\int_{\Omega}\left[ \frac{|u_i|^2}{2}-
    F_i(u_i)\right]\,dx\geq \sum_{i=1}^k\left( \frac{|A_i|^2}2-
    F_i(A_i)\right)|\{x\in\Omega:u_i(x)> 0\}|\\
  \notag& \geq -\mu_{i_0}\sum_{i=1}^k|\{x\in\Omega:u_i(x)> 0\}|\geq
  -\mu_{i_0}|\Omega|=J_{\Omega}(U_0).
\end{align}
On the other hand for any nontrivial $k$-uple
$U=(U_1,\dots,U_k)\in{\mathcal U}$ there exists $j$ such that $|\nabla
u_j|\not\equiv 0$ and hence the inequality in the first line of
\eqref{eq:1} is strict. Therefore $J_{\Omega}(U)>
J_{\Omega}(U_0)$ and $U$ cannot be a global minimizer.
\end{pfn}

A nontrivial solution to the local minimization problem will be
provided by a limit configuration of solutions to the competitive
system. To this aim we shall perform the asymptotic analysis of the
solutions to (\ref{e:nostro}) found in Theorem \ref{t:sistemi2} as
$\vk\to+\infty$.

%\subsection{Asymptotic analysis}
%We are going to prove that, if $U^{\e,\vk}$ is a solution to
% (\ref{e:nostro}) close to $W$ and total energy $c_{\e,\vk}$, all the
% species survive as $\vk\to\infty$ in a configuration of segregation,
% which is a local minimizer of the free energy $J_{\Omega_\e}$.  This
% analysis will directly provide the proof both to
\medskip\noindent
\begin{pfn}{Theorem \ref{t:partizione} and
\ref{t:sistemi}}
Let  $U^{\e,\vk}=(u_1^{\e,\vk},...,u_k^{\e,\vk})$ be the solution
of   system \eqref{e:nostro} obtained in  Theorem \ref{t:sistemi2} by minimizing
$I_{\e,\vk}$ on $B^\e_{\delta}(W)$,  hence $I_{\e,\vk}(U^{\e,\vk})=c_{\e,\vk}$
as in Theorem \ref{t:auxiliary}.
In particular
\begin{equation}\label{eq:8}
c_{\e,\vk}\geq \frac12\|U^{\e,\vk}\|_{(H^1(\Omega_\e))^k}^2-
|\Omega_\e|\sum_i \max_{t\in[0,A_i]}|F_i(t)|.
\end{equation}
For every $U\in{\mathcal U}_\e\cap B^\delta_\e(W)$, define
$\widetilde U$ by setting $\tilde u_i(x)=\min\{u_i(x), A_i\}$.
Then the following inequalities hold
$$
J_{\Omega_\e}(U)\geq J_{\Omega_\e}({\widetilde U})=I_{\e,\k}(U)\geq c_{\e,\vk},
$$
implying
\begin{equation}\label{eq:9}
\lambda_{\e}^{\delta}\geq c_{\e,\k}.
\end{equation}
From (\ref{eq:8}) and (\ref{eq:9}) we obtain that
\[
\|U^{\e,\vk}\|_{(H^1(\Omega_\e))^k}^2\leq 2
c_{\e,\vk}+2|\Omega_\e|\sum_i \max_{t\in[0,A_i]}|F_i(t)| \leq
2\lambda_{\e}^{\delta}+2|\Omega_\e|\sum_i \max_{t\in[0,A_i]}|F_i(t)|.
\]
Hence $u_i^{\e,\vk}$ is bounded in $H^1(\Omega_\e)$ uniformly with
respect to $\vk$, then there exists a weak limit $v_i^\e$ such that,
up to subsequences, $u_i^{\e,\vk}\rightharpoonup v_i^{\e}$ in
$H^1(\Omega_\e)$ as $\vk\to+\infty$.  Also, by lower semicontinuity of
the norm, we learn that $V^\e\in B_\delta^\e(W)$, hence, by
(\ref{eq:nontrivi}), $v_i^\e\not\equiv 0$ for all $i$.  Let us now
multiply the equation of $u_i^{\e,\vk}$ times $u_i^{\e,\vk}$ on
account of the boundary conditions: then
\begin{eqnarray*}
  \vk\int_{\Omega}(u^{\e,\vk}_i)^2\sum_{j\neq
    i}(u^{\e,\vk}_j)^2\quad\text{is bounded uniformly in }\vk,
\end{eqnarray*}
hence
\[
\int_{\Omega}(u^{\e,\vk}_i)^2\sum_{j\neq i}(u^{\e,\vk}_j)^2\to
0,\qquad \text{as }\vk\to\infty.
\]
By the pointwise convergence $u_i^{\e,\vk}(x)\to v_i^{\e}(x)$ a.e.
$x\in\Omega_\e$, we infer that $v_i^\e(x)\geq 0$ and
$v_i^\e(x)\cdot v_j^\e(x)=0$ for almost every $x$, hence
$V^\e\in\mathcal U_\e$.

Also, by the positivity of the interaction term, we know that
$c_{\e,\vk}\leq c_{\e,\vk'}$ when $\vk\leq \vk'$: hence the
sequence of critical levels $c_{\e,\vk}$ converges to some
$\lambda\leq \lambda_{\e}^{\delta}$ as $\vk\to +\infty$.  Since by
the Dominated Convergence Theorem (recall that $0\leq
u_i^{\e,\vk}\leq A_i$)
$$
\int_{\Omega_\e} F_i(u_i^{\e,\vk})\,dx =\int_{\Omega_\e}
\widetilde F_i(u_i^{\e,\vk})\,dx\to \int_{\Omega_\e}
 F_i(v_i^{\e})\,dx,\qquad \vk\to\infty,
$$
the following chain of inequalities holds:
\begin{align*}
  \lambda_{\e}^{\delta}&\geq\lim_{\vk\to\infty}c_{\e,\vk}
  =\lim_{\vk\to\infty}I_{\e,\vk}(U^{\e,\vk})\\
  & = \limsup_{\vk\to\infty}\bigg[
  \sum_{i=1}^k\bigg\{\frac12\|u_i^{\e,\vk}\|^2_{H^1(\Omega_\e)}
  -\int_{\Omega_\e} \widetilde F_i(u_i^{\e,\vk})\,dx\bigg\}+
  \vk\sum_{\substack{i,j=1\\i\neq j}}^k\int_{\Omega}
  (u^{\e,\vk}_i)^2(u^{\e,\vk}_j)^2\bigg]\\
  &\geq \limsup_{\vk\to\infty}
  \sum_{i=1}^k\bigg\{\frac12\|u_i^{\e,\vk}\|^2_{H^1(\Omega_\e)}
  -\int_{\Omega_\e} \widetilde F_i(u_i^{\e,\vk})\,dx\bigg\}\\
  &\geq \liminf_{\vk\to\infty}
  \sum_{i=1}^k\bigg\{\frac12\|u_i^{\e,\vk}\|^2_{H^1(\Omega_\e)}
  -\int_{\Omega_\e} \widetilde F_i(u_i^{\e,\vk})\,dx\bigg\}\\
  &\geq \sum_{i=1}^k\bigg\{\frac12\|v_i^{\e}\|^2_{H^1(\Omega_\e)}
  -\int_{\Omega_\e} F_i(v_i^{\e})\,dx\bigg\}= J_{\Omega_\e}(V^\e)\geq
  \lambda_{\e}^{\delta}.
\end{align*}
Therefore all the above inequalities are indeed equalities. In
particular $J_{\Omega_\e}(V^\e)=\lambda_{\e}^{\delta}$, meaning that $V^\e$
solves ${\boldsymbol{(P_\e)}}$ on $B_\delta^\e(W)$,
giving the proof of Theorem \ref{t:partizione}.

 Moreover
$\lim_{\vk\to+\infty}\|U^{\e,\vk}\|_{(H^1(\Omega_\e))^k}=
\|V^{\e}\|_{(H^1(\Omega_\e))^k}$ which, together with weak
convergence, implies that the convergence $U^{\e,\vk}\to V^\e$ is
actually strong in $(H^1(\Omega_\e))^k$. We also deduce that
\[
\lim_{\vk\to\infty}\vk\int_{\Omega}(u^{\e,\vk}_i)^2\sum_{j\neq
  i}(u^{\e,\vk}_j)^2=0.
\]
The proof of the Theorem \ref{t:sistemi} is thereby complete.
\end{pfn}

\subsection{Extremality conditions}\label{sec:extr-cond}
Once the existence of a solution for the optimal partition problem
${\boldsymbol{(P_\e)}}$ is known, we can appeal to \cite{ctv-var} to
derive some interesting properties of $U_\e$.  In particular, since
$U^\e$ is a local minimizer of the free energy $J_{\Omega_\e}$ we can
prove that its components are solution of a remarkable system of
differential inequalities.
\begin{Theorem}\label{t:extremality} Let
  $U^\e\in B^\delta_\e(W)$ be a solution to problem
  ${\boldsymbol{(P_\e)}}$. Then $U^\e$ is a solution of the $2k$
  distributional inequalities \eqref{2kvar}, namely, for every $i$ and
  every $\phi \in H^1(\Omega_\e)$ such that $\phi\geq 0$ a.e. in
  $\Omega_\e$, there holds
\begin{align*}
 \begin{cases}
{\displaystyle   \int_{\Omega_\e}}\big(\nabla u^\e_i\nabla\phi+u_i^\e\phi-f_i(u^\e_i)
\phi\big)\,dx\leq 0,\\[10pt]
{\displaystyle \int_{\Omega_\e}}\big(\nabla \widehat
u^\e_i\nabla\phi+\widehat u_i^\e\phi-\widehat f_i(\widehat
u^\e_i)\phi\big)\,dx\geq 0,
\end{cases}
\end{align*}
where $\widehat u_i = u_i-\sum_{h\neq i}u_h$ and $\widehat
f(\widehat u_i) =f_i(u_i)-\sum_{j\neq i} f_j(u_j)$.
\end{Theorem}
The proof can be obtained as in \cite[Theorem 5.1]{ctv-var}, the only
difference being that here we are dealing with local (and not
global) minima of the free energy.
  For the reader's convenience, we
  sketch the main steps.

 \medskip\noindent
 \begin{pf}   We argue by contradiction, hence, to prove the
  first inequality, we assume that there exists one index $j$ and
  $\phi\in H^1(\Omega_\e)$ such that $\phi\geq 0$ and
\begin{equation}\label{e:ass1}
\int_{\Omega_\e}\big(\nabla u^\e_j\nabla
\phi+u_j^\e\phi-f_j(u_j^\e)\phi\big)>0.
\end{equation}
For $t\in (0,1)$ we consider $V=(v_1,\dots,v_k)$ defined as
$$
v_{i}=\left\{\begin{array}{ll}u^\e_{i} &\mbox{ if $i\neq j$}
    \\(u_i^\e-t\phi)^+&\mbox{ if $i=j$.}\end{array}\right.
$$
We notice that $V\in{\mathcal U}_\e$. Moreover, since that map
$z\mapsto [z]^+$ is continuous from $H^1(\Omega_\e)$ to
$H^1(\Omega_\e)$ and $U^\e\in B_\delta^\e(W)$, we learn that $V\in
B_\delta^\e(W)$ for all $t$ small enough.  In light of \eqref{e:ass1}
it is immediate to check that
$J_{\Omega_\e}(V)<J_{\Omega_\e}(U^\e)=\min\{J_{\Omega_\e}(U),\,U\in
B_\delta^\e(W)\cap{\mathcal U}_\e\}$ for $t$ small enough, a
contradiction.  Let now $j$ and $\phi\in H^1(\Omega_\e)$, $\phi\geq
0$, such that
$$
\int_\Omega\big(\nabla \widehat u^\e_j\nabla\phi+u^\e_i\phi -\widehat
f(\widehat u^\e_j)\phi\big)<0.
$$
Again, we show that the value of the functional can be lessen by
replacing $U$ with an appropriate new function $V$ close to $W$.
This is defined as  $V=(v_1,\dots,v_k)$ with
$$
v_i=\left\{
\begin{array}{ll}
\left(\widehat u_j+t\phi\right)^+, &\mbox{if $i=j$}\\
\left(\widehat u_j+t\phi\right)^-\chi_{\supp(u_i)}, &\mbox{if
$i\neq j$}.
\end{array}\right.
$$
Simple computations lead to
$$
J_{\Omega_\e}(V)-J_{\Omega_\e}(U^\e)=t\int_{\Omega_\e}\big(\nabla \widehat
u^\e_j\nabla\phi+\widehat u^\e_i\phi-\widehat f_j(\widehat
u^\e_j)\phi\big)+o(t),
$$
which leads to a contradiction if $t$ is small enough.
\end{pf}

% Furthermore, if the potential energies $\frac12 u_i^2-F_i(u_i)$ are
% convex as for the logistic internal growth $f_i(s)=\lambda
% u-u|u|^{p-1}$ $(p>1)$, we can easily adapt to the present case the
% proof of \cite[Theorem 4.2]{ctv-var} to guarantee that $U_\e$ is the
% \emph{unique} minimizer of $J_{\Omega_\e}$ close to $W$.
%\begin{Theorem}\label{teo2} Assume that (F1)--(F3) hold and that
%\begin{eqnarray*}
%f_i'(s)\leq 1,\; \forall s\in\R.
%\end{eqnarray*}
%for all $i=1,\dots,k$. Then, problem ${\boldsymbol{(P_\e)}}$ has a
% unique local minimizer in $B_\e^\delta(W)$.
%\end{Theorem}
%It is worth noticing that, due to the positivity of the interaction
% term, \emph{any} solution of the system \eqref{e:nostro} with
% non-negative components do satisfies inequalities \eqref{2kvar}.
% Hence \emph{any} $U\in \mathcal U_\e$ obtained as weak limit of
% solutions to the system as $\vk\to\infty$, is also a solution of
% \eqref{2kvar}.

\section{Conclusions and final remarks}
As a final step of our study, we have proved the existence of an
element $(u_1,\dots,u_k)$ with $k$ non-trivial components in the
functional class
$${\mathcal S}(\Omega)=\left\{
\begin{array}{c}
(u_1,\cdots,u_k)\in (H^1(\Omega))^k:\
 u_i\geq0,\, u_i\not\equiv 0,\, u_i\cdot u_j=0\mbox{ if }\,i\neq j,\\[5pt]
  \int_{\Omega}\big(\nabla u_i\nabla\phi+u_i\phi-f_i(u_i)\phi\big)
\leq 0\text{ and }
  \int_{\Omega}\big(\nabla \widehat u_i\nabla\phi+\widehat
  u_i\phi-\widehat f(\widehat u_i)\phi\big)\geq 0\\[5pt]
\text{for every }i=1,\dots,k \text{ and }\phi\in H^1(\Omega)
 \text{ such that }\phi\geq0\text{ a.e. in }\Omega
 \end{array}\right\}
$$
when $\Omega=\Omega_\e$ with small $\e$.

In particular, by choosing test functions $\phi$ with compact support
in $\Omega_\e$, we learn that any element of ${\mathcal S}(\Omega_\e)$
is a solution (in distributional sense) of the following $2k$
differential inequalities:
\begin{equation}\label{e:diffineq}
\begin{cases}
  -\Delta u_i\leq f_i(u_i),\quad &\text{in }\Omega_\e,\\
  -\Delta \widehat u_i\geq \widehat f(\widehat u_i),
  &\text{in }\Omega_\e.
\end{cases}
\end{equation}
By appealing to the interior regularity theory developed in
\cite[Section 8]{ctv-var}, we know that any $u_i$ is locally Lipschitz
continuous and, in particular, the set $\omega_i=\{x\in\Omega_\e:\;
u_i(x)>0\}$ is an open (nonempty) set. Hence by \eqref{2kvar} we
obtain that $u_i\big|_{\omega_i}$ is solution of
$$-\Delta u_i+u_i=f_i(u_i),\qquad \text{in }\omega_i,$$
subject to the boundary condition $$\frac{\partial
  u_i}{\partial\nu}=0,\qquad\text{on }\partial\Omega_\e\cap \omega_i.$$
This suggest that the validity of \eqref{2kvar} not only implies the
 differential inequalities \eqref{e:diffineq} in
$\Omega_\e$, but it also contains boundary conditions on
$\partial\Omega_\e$ in some Neumann form, the major difficulty being
to give functional sense to ``$\frac{\partial u_i}{\partial\nu}$'' on the
whole of $\partial\Omega_\e$.  A rigorous analysis of this point
requires the development of a regularity theory for the class
$\mathcal S(\Omega)$ up to the boundary, that will be object of future
studies.\\

\noindent\textbf{Acknowledgements.} We are indebted to 
Prof. Susanna Terracini for her interesting comments and
suggestions.


\begin{thebibliography}{99}



\bibitem{am} A. Ambrosetti, A. Malchiodi, \emph{Nonlinear analysis and
    semilinear elliptic problems}, Cambridge Studies in Advanced
  Mathematics, 104. Cambridge University Press, Cambridge, 2007.

\bibitem{arrieta} J. M. Arrieta, A. Carvalho, G. Lozada-Cruz,
  \emph{Dynamics in dumbbell domains. I. Continuity of the set of
    equilibria}, J. Differential Equations, 231 (2006), no. 2,
  551--597.

\bibitem{cf} M. Conti, V. Felli, \emph{Coexistence and segregation for
    strongly competing species in special domains}, Interfaces Free
  Bound., 10 (2008), 173--195.

\bibitem{ctv2} M. Conti, S. Terracini, G. Verzini, \emph{An optimal
    partition problem related to nonlinear eigenvalues}, Journal of
  Funct. Anal., 198 (2003), no. 1, 160--196.

\bibitem{ctv-var} M. Conti, S. Terracini, G. Verzini, \emph{A
    variational problem for the spatial segregation of
    reaction--diffusion systems}, Indiana Univ. Math. J., 54 (2005),
  no. 3, 779--815.

\bibitem{ctv-fucik} M. Conti, S. Terracini, G. Verzini, \emph{On a
    class of optimal partition problems related to the Fu\v{c}\'\i k
    spectrum and to the monotonicity formulae}, Calculus of
  Variations, 22 (2005), no. 1, 45--72.

\bibitem{ctv-asymp} M. Conti, S. Terracini, G. Verzini, \emph{
    Asymptotic estimates for the spatial segregation of competitive
    systems}, Advances in Mathematics, 195 (2005), no. 2, 524--560.

\bibitem{ctv-uniq} M. Conti, S. Terracini, G. Verzini,
  \emph{Uniqueness and Least Energy Property for Solutions to Strongly
    Competing Systems}, Interfaces Free Bound., 8 (2006), no. 4,
  437--446.

\bibitem{d} E.N. Dancer, \emph{The effect of domain shape on the
    number of positive solutions of certain nonlinear equations}, J.
  Differential Equations, 74 (1988), 120--156.

\bibitem{dd3} E.N. Dancer, Y.H. Du, \emph{Competing species equations
    with diffusion, large interactions, and jumping nonlinearities},
  J. Differential Equations, 114 (1994), 434--475.

\bibitem{ddh} E.N. Dancer, Y.H. Du, \emph{Positive solutions for a
    three-species competition system with diffusion. {I}. {G}eneral
    existence results}, Nonlinear Anal., 24 (1995), no. 3, 337--357.

\bibitem{dd} E.N. Dancer, Y.H. Du, \emph{Positive solutions for a
    three-species competition system with diffusion. {II}. {T}he case
    of equal birth rates}, Nonlinear Anal., 24 (1995), no. 3,
  359--373.

\bibitem{dg1} E.N. Dancer, Z.M. Guo, \emph{Uniqueness and stability
    for solutions of competing species equations with large
    interactions}, Comm. Appl. Nonlinear Anal., 1 (1994), 19--45.

\bibitem{dhmp} E.N. Dancer, D. Hilhorst, M. Mimura, L.A. Peletier,
  \emph{Spatial segregation limit of a competition--diffusion system},
  European J. Appl. Math., 10 (1999), 97--115.

\bibitem{eflg} J.C. Eilbeck, J.E. Furter, J. L{\'o}pez-G{\'o}mez,
  \emph{Coexistence in the competition model with diffusion},
  J. Differential Equations, 107 (1994), no. 1, 96--139.

\bibitem{gl} C. Gui, Y. Lou, \emph{Uniqueness and nonuniqueness of
    coexistence states in the {L}otka-{V}olterra competition model},
  Comm. Pure Appl. Math., 47 (1994), no. 12, 1571--1594.

\bibitem{kw} K. Kishimoto, H.F. Weinberger, \emph{The spatial
    homogeneity of stable equilibria of some reaction--diffusion
    system on convex domains}, J. Differential Equations, 48 (1985), 15--21.

\bibitem{kl} P. Korman, A. Leung, \emph{On the existence and uniqueness
    of positive steady states in Lotka--Volterra ecological models
    with diffusion}, Appl. Anal., 26 (1987), 145--160.

\bibitem{lm} A.C. Lazer, P.J. McKenna, \emph{On steady state solutions
    of a system of reaction--diffusion equations from biology},
  Nonlinear Anal. TMA, 6 (1982), 523--530.

\bibitem{mhmm} {H. Matano, M. Mimura}, \emph{Pattern formation in
    competition-diffusion systems in nonconvex domains}, Publ. Res.
  Inst. Math. Sci., 19 (1983), no. 3, 1049--1079.

\bibitem{skt} {N. Shigesada, K. Kawasaki, E. Teramoto}, \emph{The
    effects of interference competition on stability, structure and
    invasion of a multispecies system}, J. Math. Biol., 21 (1984),
  no. 2, 97--113.

%  \bibitem{sweers} {G. Sweers}, \emph{A sign-changing global minimizer
%      on a convex domain}, Progress in partial differential equations:
%    elliptic and parabolic problems (Pont-\`a-Mousson, 1991), 251--258,
%    Pitman Res. Notes Math. Ser., 266, Longman Sci. Tech., Harlow,
%    1992.



\end{thebibliography}
\end{document}